\documentclass[11pt]{article}
\usepackage{graphicx}
\usepackage{amssymb}
\usepackage{color}

\newcommand{\Bx}{\mathbf{x}}
\newcommand{\Ox}{\Omega_\Bx}

\newcommand{\p}{\partial}

\newcommand{\myref}[1]{(\ref{#1})}
\newcommand{\hs}[2]{\left\|#1\right\|_{H^{#2}(\Omega)}}

\newtheorem{theorem}{\textbf{Theorem}}[section]
\newtheorem{lemma}{\textbf{Lemma}}[section]
\newtheorem{remark}{\textbf{Remark}}[section]

\newtheorem{definition}{\textbf{Definition}}[section]

\newenvironment{proof}{\noindent\textbf{Proof\ }}{\hspace*{\fill}$\Box$\medskip}
\begin{document}
\title{On the Local Well-posedness of a 3D Model for Incompressible Navier-Stokes Equations with Partial Viscosity}
\author{Thomas Y. Hou\thanks{Applied and Comput. Math, Caltech,
Pasadena, CA 91125. {\it Email: hou@acm.caltech.edu.}} \and
Zuoqiang Shi\thanks{Applied and Comput. Math, Caltech,
Pasadena, CA 91125. {\it Email: shi@acm.caltech.edu.}} \and
Shu Wang\thanks{College of Applied Sciences, Beijing University of
Technology, Beijing 100124, China.
{\it Email: wangshu@bjut.edu.cn}}
}
\date{\today}
\maketitle

\begin{abstract}
  In this short note, we study the local well-posedness of a 3D model for incompressible Navier-Stokes equations with partial viscosity. This model
was originally proposed by Hou-Lei in \cite{HouLei09a}. In a recent paper,
we prove that this 3D model with partial viscosity will develop a finite
time singularity for a class of initial condition using a mixed
Dirichlet Robin boundary condition. The local well-posedness analysis
of this initial boundary value problem is more subtle than the
corresponding well-posedness analysis using a standard boundary
condition because the Robin boundary condition we consider is
non-dissipative. We establish the local well-posedness of this
initial boundary value problem by designing a Picard iteration
in a Banach space and proving the convergence of the Picard
iteration by studying the well-posedness property of the heat
equation with the same Dirichlet Robin boundary condition.
\end{abstract}

\section{Introduction}
In this short note, we prove the local well-posedness of the 3D model with
partial viscosity. The 3D model with partial viscosity has
the following form:
 \begin{eqnarray}
\label{model-partial-vis}
&&\left\{\begin{array}{rcl}
u_t &=&2u\psi_z\\
\omega_t&=&(u^2)_z+\nu \Delta \omega\\
-\Delta \psi&=&\omega\\
\end{array}\right., \quad (\Bx,z)\in \Omega=\Ox\times(0,\infty),
\end{eqnarray}
where $\Ox = (0,a)\times (0,a)$.
Let $\Gamma =\{ (\Bx,z) \; | \; \Bx \in \Ox, \; z=0\}$.
The initial and boundary conditions for \myref{model-partial-vis}
are given as following:
\begin{eqnarray}
\label{BC-w-vis-1}
&&\omega|_{\p\Omega\backslash\Gamma}=0,\quad \left(\omega_z+\gamma \omega\right)|_\Gamma=0, \\
\label{BC-psi-vis-1}
&&\psi|_{\p\Omega\backslash\Gamma}=0,\quad \left(\psi_z+\beta \psi\right)|_\Gamma=0, \\
&& \omega|_{t=0}=\omega_0(\Bx,z),\quad u|_{t=0}=u_0(\Bx,z).
\label{IC-psi-w-1}
\end{eqnarray}

This 3D model with viscosity in both $u$ and $\omega$ components
was first proposed by Hou and Lei in \cite{HouLei09a}.
The only difference between this 3D model and the reformulated
Navier-Stokes equations is that convection term is neglected
in the model. If one adds the convection term back to the 3D model,
one would recover the full Navier-Stokes equations. This model
preserves almost all the properties of the full 3D Navier-Stokes
equations. Despite the striking similarity at the theoretical level
between the 3D model and the Navier-Stokes equations, the former seems
to have a very different behavior from the full Navier-Stokes equations.
In a recent paper \cite{HSW11}, we prove that
the above 3D model with partial viscosity develops a finite time
singularity for a class of initial condition using a mixed Dirichlet
Robin boundary condition.

The analysis of finite time singularity formation of the 3D model
\cite{HSW11} uses the local well-posedness result of the 3D model.
The local well-posedness of the 3D model can be proved by
using a standard energy estimate and a mollifier if there is
no boundary or if the boundary condition is a standard one,
see e.g. \cite{MB02}. For the mixed Dirichlet Robin
boundary condition we consider here, the analysis is a bit more
complicated since the mixed Dirichlet Robin condition gives
rise to a growing eigenmode.

There are two key ingredients in our local well-posedness analysis.
The first one is to design a Picard iteration for the 3D model.
The second one is to show that the mapping that generates the
Picard iteration is a contraction mapping and the Picard
iteration converges to a fixed point of the Picard mapping
by using the Contraction Mapping Theorem. To establish the contraction
property of the Picard mapping, we need to use the well-posedness
property of the heat equation with the same Dirichlet Robin
boundary condition as $\omega$. The well-posedness analysis
of the heat equation with a mixed Dirichelet Robin boundary
has been studied in the literature. The case of $\gamma >0$ is
more subtle because there is a growing eigenmode. Nonetheless,
we prove that all the essential regularity properties of the
heat equation are still valid for the mixed Dirichlet Robin
boundary condition with $\gamma>0$.

\section{The main result}
The local existence result of our 3D model with partial viscosity is
stated in the following theorem.
\begin{theorem}
  \label{local-exist-vis}
Assume that $u_0 \in H^{s+1}(\Omega)$, $\omega_0 \in
H^{s}(\Omega)$ for some $s>3/2$, $u_0|_{\partial
\Omega}=u_{0z}|_{\partial \Omega}= 0 $ and $\omega_0$ satisfies
\myref{BC-w-vis-1}. Moreover, we assume that $\beta \in S_\infty$
(or $S_b$) as defined in Lemma \ref{reisz}. Then there exists a
finite time $T=T\left(\|u_0\|_{H^{s+1}(\Omega)},
\|\omega_{0}\|_{H^{s}(\Omega)}\right)>0$ such that the system
\myref{model-partial-vis} with boundary condition
\myref{BC-w-vis-1},\myref{BC-psi-vis-1} and initial data
\myref{IC-psi-w-1} has a unique solution, $u\in
C([0,T],H^{s+1}(\Omega))$, $\omega\in C([0,T],H^{s}(\Omega))$ and
$\psi\in C([0,T],H^{s+2}(\Omega))$.
\end{theorem}
The local well-posedness analysis relies on the following local
well-posedness of the heat equation and the elliptic equation with
mixed Dirichlet and Robin boundary conditions. First, the local
well-posedness of the elliptic equation with the mixed Dirichlet
and Robin boundary condition is given by the following lemma
\cite{HSW11}:
\begin{lemma}
\label{reisz}
There exists a unique solution $v\in H^{s}(\Omega)$
to the boundary value problem:
\begin{eqnarray}
\label{eqn-Laplace}
&&-\Delta v = f,\quad (\Bx,z)\in \Omega, \\
\label{mix bc}
&&v|_{\partial \Omega\backslash \Gamma} =0,
\quad (v_z+\beta v)|_\Gamma =0,
\end{eqnarray}
if $\beta \in S_\infty \equiv \{\beta \;|\;
\beta\ne \frac{\pi|k|}{a} \;\;\mbox{for all} \; k\in \mathbb{Z}^2\}$,
$f\in H^{s-2}(\Omega)$ with $s \ge 2$ and
$f|_{\p\Omega\backslash\Gamma}=0$. Moreover we have
\begin{eqnarray}
\|v\|_{H^s(\Omega)}\le C_s \|f\|_{H^{s-2}(\Omega)},
\end{eqnarray}
where $C_s$ is a constant depending on $s$, $|k|=\sqrt{k_1^2+k_2^2}$.
\end{lemma}

\begin{definition}
\label{def-k}
Let $\mathcal{K}: H^{s-2}(\Omega)\rightarrow H^{s}(\Omega)$ be a linear operator defined as following:
\begin{eqnarray*}
  \mbox{for all}\quad f\in H^{s-2}(\Omega),\quad \mathcal{K}(f) \;\mbox{is the solution of the boundary value problem
\myref{eqn-Laplace}-\myref{mix bc}}.
\end{eqnarray*}
\end{definition}
It follows from Lemma \ref{reisz} that for any $f\in H^{s-2}(\Omega)$,
we have
\begin{eqnarray}
\label{bound k}
  \|\mathcal{K}(f)\|_{H^s(\Omega)}\le  C_s\|f\|_{H^{s-2}(\Omega)} .
\end{eqnarray}

For the heat equation with the mixed Dirichlet and Robin boundary
condition, we have the following result.
\begin{lemma}
\label{heat-kernel}
There exists a unique solution $\omega\in C([0,T];H^{s}(\Omega))$
to the initial boundary value problem:
\begin{eqnarray}
\label{eqn-heat}
&&\omega_t=\nu \Delta \omega,\quad (\Bx,z)\in \Omega, \\
&&\omega|_{\partial \Omega\backslash \Gamma} =0, \quad
(\omega_z+\gamma \omega)|_\Gamma =0,\label{BC-robin}\\
&&\omega|_{t=0}=\omega_0(\Bx,z).\label{ic-omega}
\end{eqnarray}
for $\omega_0\in H^{s}(\Omega)$ with $s > 3/2$. Moreover we have
the following estimates in the case of $\gamma > 0$
\begin{eqnarray}
\label{heat-estimate}
\|\omega (t)\|_{H^s(\Omega)}\le
C(\gamma,s)e^{\nu\gamma^2 t} \|\omega_0\|_{H^{s}(\Omega)}, \quad
t \geq 0,
\end{eqnarray}
and
\begin{eqnarray}
\label{heat-estimate-1}
\|\omega (t)\|_{H^s(\Omega)}\le
C(\gamma,s,t) \|\omega_0\|_{L^2(\Omega)}, \quad
t >  0.
\end{eqnarray}
\end{lemma}
\begin{remark}
We remark that the growth factor $e^{\nu \gamma^2 t}$ in
\myref{heat-estimate} is
absent in the case of $\gamma \leq 0$ since there is no growing
eigenmode in this case.
\end{remark}

\begin{proof}
First, we prove the solution of the system
\myref{eqn-heat}-\myref{ic-omega} is unique.  Let
$\omega_1, \; \omega_2 \in H^s(\Omega)$ be two smooth solutions
of the heat equation for $0 \le t < T$ satisfying the same initial
condition and the Dirichlet Robin boundary condition.
Let $\omega= \omega_1 - \omega_2$. We will prove that
$\omega = 0$ by using an energy estimate and
the Robin boundary condition at $\Gamma$:
\begin{eqnarray}
  \label{unique-heat}
  \frac{1}{2}\frac{d}{dt}\int_\Omega \omega^2 d\Bx dz&=&\nu\int_\Omega \omega \Delta \omega d\Bx dz\nonumber\\
&=&-\nu\int_\Omega |\nabla \omega|^2 d\Bx dz-\nu\int_{\Gamma} \omega \omega_z d\Bx\nonumber\\
&=&-\nu\int_\Omega |\nabla \omega|^2 d\Bx dz+\nu\gamma\int_{\Gamma} \omega^2 d\Bx\nonumber\\
&=&-\nu\int_\Omega |\nabla \omega|^2 d\Bx dz-\nu\gamma\int_{\Gamma}\int_z^\infty\left(\omega^2\right)_z dz d\Bx\nonumber\\
&=&-\nu\int_\Omega |\nabla \omega|^2 d\Bx dz-2\nu\gamma\int_{\Gamma}\int_z^\infty\omega\omega_z d\Bx dz\nonumber\\
&\le & -\nu\int_\Omega |\nabla \omega|^2 d\Bx dz+
\frac{\nu}{2}\int_{\Omega}|\omega_z|^2 d\Bx dz+ 2\nu \gamma^2
\int_{\Omega}\omega^2 d\Bx dz \nonumber\\
&\le & -\frac{\nu}{2}\int_\Omega |\nabla \omega|^2 d\Bx dz+ 2\nu \gamma^2
\int_{\Omega}\omega^2 d\Bx dz,
\end{eqnarray}
where we have used the fact that the smooth solution of the heat 
equation $\omega$ decays to zero as $z \rightarrow \infty$.
Thus, we get
\begin{eqnarray}
  \frac{1}{2}\frac{d}{dt}\int_\Omega \omega^2 d\Bx dz \le 2\nu \gamma^2
\int_{\Omega}\omega^2 d\Bx dz.
\end{eqnarray}
It follows from Gronwall's inequality
\begin{eqnarray}
  e^{-4\nu\gamma^2 t}\int_\Omega \omega^2 d\Bx dz\le
  \int_\Omega \omega_0^2 d\Bx dz  = 0,
\end{eqnarray}
since $\omega_0 = 0$. Since $\omega \in H^s(\Omega)$ with
$s>3/2$, this implies that $\omega = 0$  for $0 \le t < T$
which proves the uniqueness of smooth solutions for the
heat equation with the mixed Dirichlet Robin boundary condition.

Next, we will prove the existence of the solution by constructing a solution explicitly.
Let $\eta(\Bx,z,t)$ be the solution of the following initial boundary
value problem:
\begin{eqnarray}
\label{eqn-heat-unbound}
&&\eta_t=\nu \Delta \eta,\quad (\Bx,z)\in \Omega, \\
&&\eta|_{\partial \Omega} =0, \quad \eta|_{t=0}=\eta_0(\Bx,z),
\end{eqnarray}
and let $\xi(\Bx,t)$ be the solution of the following PDE in $\Omega_\Bx$:
\begin{eqnarray}
\label{eqn-heat-unbound-xi}
&&\xi_t=\nu \Delta_\Bx \xi+\nu\gamma^2\xi,\quad \Bx\in \Omega_\Bx, \\
&&\xi|_{\partial \Omega_\Bx} =0, \quad \xi|_{t=0}=\overline{\omega}_0(\Bx),
\end{eqnarray}
where $\Delta_\Bx=\frac{\p^2}{\p x_1^2}+\frac{\p^2}{\p x_2^2}$ and
$ \overline{\omega}_0(\Bx)=2\gamma\int_0^\infty \omega_0(\Bx,z)e^{-\gamma z}dz$.
From the standard theory of the heat equation, we know that $\eta$ and $\xi$
both exist globally in time.

We are interested in the case when the initial value
$\eta_0(\Bx,z)$ is related to $\omega_0$ by solving the
following ODE as a function of $z$ with ${\Bx}$ being fixed as
a parameter:
\begin{eqnarray}
  \label{eq-eta0}
 -\frac{1}{\gamma}\eta_{0z}+\eta_0=\omega_0(\Bx,z)- \overline{\omega}_0(\Bx)e^{-\gamma z},\quad \eta_0(\Bx,0)=0.
\end{eqnarray}
Define
\begin{eqnarray}
  \label{sol-half}
  \omega(\Bx,z,t)\equiv-\frac{1}{\gamma}\eta_z+\eta+\xi(\Bx,t)e^{-\gamma z},\quad (\Bx,z)\in \Omega .
\end{eqnarray}
It is easy to check that
$\omega$ satisfies the heat equation for $t >0$ and the
initial condition. Obviously, $\omega$ also satisfies the
boundary condition on $\partial \Omega\backslash \Gamma$. To verify
the boundary condition on $\Gamma$, we observe by a direct calculation that
$(\omega_z + \gamma \omega ) |_{\Gamma}=
-\frac{1}{\gamma} \eta_{zz}|_\Gamma.$
Since $\eta(\Bx,z)|_\Gamma = 0$, we obtain by using
$\eta_t = \nu \Delta \eta$ and taking the limit as $z \rightarrow 0+$
that $\Delta \eta|_\Gamma = 0$, which implies that
$\eta_{zz}|_\Gamma = 0$. Therefore, $\omega$ also satisfies
the Dirichlet Robin boundary condition at $\Gamma$.
This shows that $\omega$ is a solution of the system
\myref{eqn-heat}-\myref{ic-omega}. By the uniqueness result that
we proved earlier, the solution of the heat equation must be given by
\myref{sol-half}.

Since $\eta$ and $\xi$ are solutions of the heat equation with
a standard Dirichlet boundary condition,
the classical theory of the heat equation \cite{Evans98} gives
the following regularity estimates:
\begin{eqnarray}
  \label{hs-heat-standard}
  \hs{\eta}{s}\le C \hs{\eta_0}{s},\quad
\left\|\xi(\Bx)\right\|_{H^s(\Omega_\Bx)}\le Ce^{\nu\gamma^2 t}\left\|\overline{\omega}_0(\Bx)\right\|_{H^s(\Omega_\Bx)}.\quad
\end{eqnarray}
Recall that $\eta_{zz}|_\Gamma=0$. Therefore, $\eta_z$ also solves the
heat equation with the same Dirichlet Robin boudary condition:
\begin{eqnarray}
\label{eqn-heat-etaz}
&&\left(\eta_z\right)_{t}=\nu \Delta \eta_z,\quad (\Bx,z)\in \Omega, \\
&&\left(\eta_z\right)_z|_{\Gamma} =0, \quad (\eta_z)|_{\p \Omega\backslash \Gamma}=0,\quad (\eta_z)|_{t=0}=\eta_{0z}(\Bx,z),
\end{eqnarray}
which implies that
\begin{eqnarray}
  \label{hs-heat-standard-etaz}
  \hs{\eta_z}{s}\le C \hs{\eta_{0z}}{s}.
\end{eqnarray}
Putting all the above estimates for $\eta$, $\eta_z$ and $\xi$ together
and using \myref{sol-half}, we obtain the following estimate:
\begin{eqnarray}
  \label{hs-heat-1}
  \hs{\omega}{s}&=&\hs{-\frac{1}{\gamma}\eta_z+\eta+\xi(\Bx,t)e^{-\gamma z}}{s}\nonumber\\
&\le&\frac{1}{\gamma}\hs{\eta_z}{s}+\hs{\eta}{s}+\hs{\xi(\Bx,t)e^{-\gamma z}}{s}\nonumber\\
&\le&C(\gamma,s)\left(\hs{\eta_{0z}}{s}+\hs{\eta_0}{s}+e^{\nu\gamma^2 t}\left\|\overline{\omega}_0(\Bx)\right\|_{H^s(\Omega_\Bx)}\right).
\end{eqnarray}

It remains to bound $\hs{\eta_{0z}}{s}$, $\hs{\eta_0}{s}$ and $\left\|\overline{\omega}_0(\Bx)\right\|_{H^s(\Omega_\Bx)}$ in terms of
 $\hs{\omega_{0}}{s}$.
By solving the ODE \myref{eq-eta0} directly, we can express
$\eta$ in terms of $\omega_0$ explicitly
\begin{eqnarray}
  \label{sol-eta0}
  \eta_0(\Bx,z)=-\gamma e^{\gamma z}\int_0^z e^{-\gamma z'}f(\Bx,z')dz'
= \gamma \int_z^\infty e^{-\gamma (z'-z)}f(\Bx,z')dz' ,
\end{eqnarray}
where $f(\Bx,z)=\omega_0(\Bx,z)-\overline{\omega}_0(\Bx)e^{-\gamma z}$ and we have used the property that
\begin{eqnarray*}
  \int_0^\infty f(\Bx,z)e^{-\gamma z}dz=0.
\end{eqnarray*}
By using integration by parts, we have
\begin{eqnarray}
  \eta_{0z}(\Bx,z)= -\gamma f(\Bx,z)+ \gamma^2 \int_z^\infty e^{-\gamma (z'-z)}f(\Bx,z')dz' =\gamma \int_z^\infty e^{-\gamma (z'-z)}f_{z'}(\Bx,z')dz'.
\end{eqnarray}
By induction we can show that for any
$\alpha=(\alpha_1,\alpha_2,\alpha_3)\ge 0$
\begin{eqnarray}
  D^{\alpha}\eta_0= \gamma \int_z^\infty e^{-\gamma (z'-z)}D^\alpha f(\Bx,z')dz'.
\end{eqnarray}
Let $K(z)=\gamma e^{-\gamma z}\chi(z)$ and $\chi(z)$ be the characteristic function
\begin{eqnarray}
\chi(z)=\left\{\begin{array}{ll}
0,& z\le0,\\
1,&z> 0.
\end{array}\right.
\end{eqnarray}
Then $D^\alpha\eta_0$ can be written in the following convolution form:
\begin{eqnarray}
  \label{eta0-conv}
  D^\alpha\eta_0(\Bx,z)=  \int_0^\infty K(z'-z)D^\alpha f(\Bx,z')dz'.
\end{eqnarray}
Using Young's inequality (see e.g. page 232 of \cite{Foland84}), we obtain:
\begin{eqnarray}
  \label{eta0-l2}
  \|D^\alpha\eta_0\|_{L^2(\Omega)}&\le& \|K(z)\|_{L^1\left(\mathbb{R}^+\right)}\|D^\alpha f\|_{L^2(\Omega)}\nonumber\\
&\le& C(\gamma) \left\|D^\alpha \omega_0-\left(-\gamma\right)^{\alpha_3}e^{-\gamma z}D^{(\alpha_1,\alpha_2)}\overline{\omega}_0(\Bx)\right\|_{L^2(\Omega)}
\nonumber\\
&\le& C(\gamma,\alpha)\left(\left\|D^\alpha \omega_0\right\|_{L^2(\Omega)}+\left\|D^{(\alpha_1,\alpha_2)}\overline{\omega}_0(\Bx)\right\|_{L^2(\Omega_\Bx)}
\right) .
\end{eqnarray}
Moreover, we obtain by using the H\"{o}lder inequality that
\begin{eqnarray}
\label{omega-l2}
  \left\|D^{(\alpha_1,\alpha_2)}\overline{\omega}_0(\Bx)\right\|_{L^2(\Omega_\Bx)}&=&
 \left(\int_{\Omega_\Bx}\left(\int_0^\infty e^{-\gamma z} D^{(\alpha_1,\alpha_2)}\omega_0(\Bx,z)dz\right)^2 d\Bx\right)^{1/2}\nonumber\\
&\le &\left(\frac{1}{2\gamma} \int_{\Omega_\Bx}\int_0^\infty\left( D^{(\alpha_1,\alpha_2)}\omega_0(\Bx,z)\right)^2 dz d\Bx\right)^{1/2}\nonumber\\
&=&  \frac{1}{\sqrt{2\gamma}} \left\|D^{(\alpha_1,\alpha_2)}\omega_0(\Bx,z)\right\|_{L^2(\Omega)} .
\end{eqnarray}
Substituting \myref{omega-l2} to \myref{eta0-l2} yields
\begin{eqnarray}
  \|D^\alpha\eta_0\|_{L^2(\Omega)}&\le&
 C(\gamma,\alpha)\left(\left\|D^\alpha \omega_0\right\|_{L^2(\Omega)}+\left\|D^{(\alpha_1,\alpha_2)}\omega_0\right\|_{L^2(\Omega)}
\right) ,
\end{eqnarray}
which implies that
\begin{eqnarray}
  \label{eta0-hs}
\hs{\eta_0}{s}\le C(\gamma,s)\hs{\omega_0}{s},\quad \forall \; s\ge 0.
\end{eqnarray}
It follows from \myref{omega-l2} that
\begin{eqnarray}
  \label{omega0-hs}
  \left\|\overline{\omega}_0(\Bx)\right\|_{H^s(\Omega_\Bx)}\le C(\gamma)\left\|\omega_0\right\|_{H^s(\Omega)},\quad \forall\; s\ge 0.
\end{eqnarray}
On the other hand, we obtain from the equation for $\eta_0$ \myref{eq-eta0}
that
\begin{eqnarray}
  \label{eta0z-hs}
  \hs{\eta_{0z}}{s}= \gamma\hs{f+\eta_0}{s}\le C(\gamma,s)\hs{\omega_0}{s},\quad \forall \; s\ge 0.
\end{eqnarray}
Upon substituting \myref{eta0-hs}-\myref{eta0z-hs} to \myref{hs-heat-1},
we obtain
\begin{eqnarray}
  \label{hs-heat}
  \hs{\omega}{s}\le C(\gamma,s)e^{\nu\gamma^2 t}\hs{\omega_0}{s},
\end{eqnarray}
where $C(\gamma,s)$ is a constant depending on $\gamma$ and $s$ only.
This proves \myref{heat-estimate}.

To prove \myref{heat-estimate-1}, we use the classical regularity
result for the heat equation with the homogeneous Dirichlet boundary
condition to obtain the following estimates for $t>0$:
\begin{eqnarray}
  \label{l2-heat-standard1}
&&  \hs{\eta}{s}\le C(t) \|\eta_0\|_{L^2(\Omega)},\\
&& \hs{\eta_z}{s}\le C(s,t) \|\eta_{0z}\|_{L^2(\Omega)},\\
  \label{l2-heat-standard3}
&&\left\|\overline{\omega}(\Bx)\right\|_{H^s(\Omega_\Bx)}\le C(s,t)e^{\nu\gamma^2 t}\left\|\overline{\omega}_0(\Bx)\right\|_{
L^2(\Omega_\Bx)},
\end{eqnarray}
where $C(s, t)$ is a constant depending on $s$ and $t$.
By combining \myref{l2-heat-standard1}-\myref{l2-heat-standard3} 
with estimates \myref{eta0-hs}-\myref{eta0z-hs}, we obtain for 
any $t>0$ that
\begin{eqnarray}
  \label{l2-heat}
  \hs{\omega}{s}&\le&C(\gamma,s,t)\left(\|\eta_{0z}\|_{L^2(\Omega)}+\|\eta_0\|_{L^2(\Omega)}
+e^{\nu\gamma^2 t}\left\|\overline{\omega}_0(\Bx)\right\|_{L^2(\Omega_\Bx)}\right)\nonumber\\
&\le& C(\gamma,s,t)\|\omega_0\|_{L^2(\Omega)},
\end{eqnarray}
where $C(\gamma,s,t)<\infty$ is a constant depending on $\gamma$, $s$ and $t$. This proves \myref{heat-estimate-1} and completes the proof of the
Lemma.
\end{proof}

We also need the following well-known Sobolev inequality \cite{Foland95}.
\begin{lemma}
\label{lemma hs}
Let $u,v \in H^s(\Omega)$ with
$s>3/2$. We have
\begin{eqnarray}
\|uv\|_{H^s(\Omega)}\le c \|u\|_{H^s(\Omega)}\|v\|_{H^s(\Omega)}.
\end{eqnarray}
\end{lemma}

Now, we are ready to give the proof of Theorem \ref{local-exist-vis}.

\begin{proof} {\bf of Theorem \ref{local-exist-vis}}
Let $v=u^2$.
First, using the definition of the operator $\mathcal{K}$
(see Definition \ref{def-k}), we can rewrite the 3D model
with partial viscosity in the following equivalent form:
 \begin{eqnarray}
\label{regs1} &&\left\{\begin{array}{rcl}
v_t &=&4v\mathcal{K}(\omega)_z\\
\omega_t&=&v_z+\nu \Delta \omega\\
\end{array}\right., \quad (\Bx,z)\in \Omega=\Ox\times(0,\infty),
\end{eqnarray}
 with the initial and boundary conditions given as follows:
\begin{eqnarray}
\label{BC-w-vis-2}
&&\omega|_{\p\Omega\backslash\Gamma}=0,\quad \left(\omega_z+\gamma \omega\right)|_\Gamma=0, \\
&& \omega|_{t=0}=\omega_0(\Bx,z)\in W^s,\quad
v|_{t=0}=v_0(\Bx,z)\in V^{s+1} , \label{IC-psi-w}
\end{eqnarray}
where $V^{s+1}=\{v\in H^{s+1}:v|_{\p\Omega}=0, v_z|_{\p\Omega}=0, v_{zz}|_{\p\Omega}=0\}$
and $W^s=\{w\in H^s:w|_{\p\Omega\backslash\Gamma}=0,
(w_z+\gamma w)|_\Gamma=0\}$.

We note that
the condition $u_0|_{\p\Omega}=u_{0z}|_{\p \Omega}=0$ implies that
$v_0|_{\p\Omega}=v_{0z}|_{\p \Omega}= v_{0zz}|_{\p \Omega}=0$ by
using the relation $v_0=u_0^2$. Thus we have $v_0\in V^{s+1}$.
It is easy to show by using the $u$-equation that the property
$u_0|_{\p\Omega}=u_{0z}|_{\p \Omega}=0$ is preserved dynamically.
Thus we have $v \in V^{s+1}$.

Define $U = (U_1, U_2)= (v,\omega)$ and
$X=C([0, T];  V^{s+1})\times C([0,T]; W^s)$ with the norm
\begin{eqnarray}
  \|U\|_{X}=\sup_{t\in [0,
T]}\hs{U_1}{s+1}+\sup_{t\in [0, T]}\hs{U_2}{s},\quad \forall U\in
X\nonumber
\end{eqnarray}
and let $S=\{U\in X:\|U\|_X\le M\}$.

Now, define the map $\Phi: X\rightarrow X$ in the following way:
let $\Phi(\tilde{v},\tilde{\omega})=(v,\omega)$, then for any $t\in [0,T]$,
\begin{eqnarray}
  \label{picard-iteration}
  v(\Bx,z,t)&=&v_0(\Bx,z,t)+4\int_0^t\tilde{v}(\Bx,z,t')\mathcal{K}(\tilde{\omega})_z(\Bx,z,t')dt' ,\\
\omega(\Bx,z,t)&=& \mathcal{L}(\tilde{v}_z,\omega_0;\Bx,z,t) ,
\end{eqnarray}
where $\omega(\Bx,z,t)=\mathcal{L}(\tilde{v}_z,\omega_0;\Bx,z,t)$ is the solution of the following equation:
  \begin{eqnarray}
\omega_t=\tilde{v}_z+\nu \Delta \omega, \quad (\Bx,z)\in \Omega=\Ox\times(0,\infty),
\end{eqnarray}
 with the initial and boundary conditions:
\begin{eqnarray}
\omega|_{\p\Omega\backslash\Gamma}=0,\quad \left(\omega_z+\gamma \omega\right)|_\Gamma=0, \quad \omega|_{t=0}=\omega_0(\Bx,z).\nonumber
\end{eqnarray}
We use the map $\Phi$ to define a Picard iteration:
$U^{k+1}=\Phi(U^k)$ with $U^0=(v_0,\omega_0)$.
In the following, we will prove that there exist $T>0$ and $M>0$
such that
\begin{itemize}
\item[1.] $U^{k}\in S,\quad$ for all $k$.
\item[2.] $\left\|U^{k+1}-U^{k}\right\|_X\le \frac{1}{2}\left\|U^{k}-U^{k-1}\right\|_X$, for all $k$.
\end{itemize}
Then by the contraction mapping theorem, there exists $U=(v,\omega)\in S$ such that $\Phi(U)=U$ which implies that
$U$ is a local solution of the system \myref{regs1} in $X$.

First, by Duhamel's principle, we have for any $g\in C([0, T]; V^{s})$ that
\begin{eqnarray}
  \label{duhamle}
  \mathcal{L}(g,\omega_0;\Bx,z,t)=\mathcal{P}(\omega_0;0,t)+\int_0^t\mathcal{P}(g;t',t)dt' ,
\end{eqnarray}
where $\mathcal{P}(g;t',t)=\tilde{g}(\Bx,z,t)$ is defined as
the solution of the following initial boundary value problem at time
$t$:
\begin{eqnarray}
  \label{heat-eq-1}
  \tilde{g}_t=\nu \Delta \tilde{g},  \quad (\Bx,z)\in \Omega=\Ox\times(0,\infty),
\end{eqnarray}
with the initial and boundary conditions:
\begin{eqnarray}
\label{BC-heat-1}
\tilde{g}|_{\p\Omega\backslash\Gamma}=0,\quad \left(\tilde{g}_z+\gamma \tilde{g}\right)|_\Gamma=0,\quad \tilde{g}(\Bx,z,t')=g(\Bx,z,t'). \label{IC-heat}
\end{eqnarray}
We observe that $g(\Bx,z,t')$ also satisfies the same boundary condition
as $\omega$ for any $0\le t'\le t$ since $g=v^k_z$ and $v^k\in V^{s+1}$.

Now we can apply Lemma \ref{heat-kernel} to conclude that
for any $t'<T$ and $t\in [t',T]$ we have
\begin{eqnarray}
  \label{Hs-p}
  \|\mathcal{P}(g;t',t)\|_{H^s(\Omega)}\le C(\gamma,s)e^{\nu\gamma^2 (t-t')}\|g(\Bx,z,t')\|_{H^s(\Omega)}.
\end{eqnarray}
which implies the following estimate for $\mathcal{L}$: for all $t\in [0,T]$,
\begin{eqnarray}
  \label{Hs-L}
 \left\|\mathcal{L}(g,\omega_0;\Bx,z,t)\right\|_{H^s(\Omega)}\le C(\gamma,s)e^{\nu\gamma^2 t}\left(\|\omega_0\|_{H^s(\Omega)}+
t\sup_{t'\in [0,t]}\|g(\Bx,z,t')\|_{H^s(\Omega)}\right).
\end{eqnarray}

Further, by
using Lemma \ref{reisz} and the above estimate \myref{Hs-L} for
the sequence $U^{k}=(v^k,\omega^k)$, we get the following estimate:
\begin{eqnarray}
  \hs{v^{k+1}}{s+1}&\le& \hs{v_0}{s+1}
  +4T\sup_{t\in [0,T]}\hs{v^k(\Bx,z,t)}{s+1}\sup_{t\in [0,T]}\hs{\mathcal{K}(\omega^k)_z(\Bx,z,t)}{s+1},\nonumber\\
&\le & \hs{v_0}{s+1}+4T\sup_{t\in [0,T]}\hs{v^k(\Bx,z,t)}{s+1}\sup_{t\in [0,T]}\hs{\omega^k(\Bx,z,t)}{s},\forall t\in [0,T]\label{estimate-vk}
\quad\quad\\
\hs{\omega^{k+1}}{s}&\le& C(\gamma,s)e^{\nu\gamma^2 t}\left(\hs{\omega_0}{s}+t\sup_{t'\in [0,t]}\hs{v^{k}_z(\Bx,z,t')}{s}\right)\nonumber\\
&\le & C(\gamma,s)e^{\nu\gamma^2 T}\left(\hs{\omega_0}{s}+T\sup_{t\in [0,T]}\hs{v^{k}}{s+1}\right),\quad \forall t\in [0,T].\quad\quad\quad
\label{estimate-wk}
\end{eqnarray}
Next, we will use mathematical induction to prove that if $T$ satisfies the
following inequality:
\begin{eqnarray}
\label{T-constrain}
   8C(\gamma,s)Te^{\nu\gamma^2 T}\left(\hs{\omega_0}{s}+2T\hs{v_0}{s+1}\right)\le 1
\end{eqnarray}
then for all $k\ge 0$ and $t\in [0,T]$, we have that
\begin{eqnarray}
  \label{hs-iteration-v}
 && \hs{v^{k}}{s+1}\le 2\hs{v_0}{s+1},\\
  \label{hs-iteration-w}
 &&\hs{\omega^{k}}{s}\le C(\gamma,s)e^{\nu\gamma^2 T}\left(\hs{\omega_0}{s}+2T\hs{v_0}{s+1}\right).
\end{eqnarray}
First of all, $U^0=(v_0,\omega_0)$ satisfies \myref{hs-iteration-v} and \myref{hs-iteration-w}. Assume $U^k=(v^k,\omega^k)$ has
this property, then for $U^{k+1}=(v^{k+1},\omega^{k+1})$, using \myref{estimate-vk} and \myref{estimate-wk}, we have
\begin{eqnarray}
    \hs{v^{k+1}}{s+1}
&\le & \hs{v_0}{s+1}+4T\sup_{t\in [0,T]}\hs{v^k(\Bx,z,t)}{s+1}\sup_{t\in [0,T]}\hs{\omega^k(\Bx,z,t)}{s}\nonumber\\
&\le & \hs{v_0}{s+1}\left(1+8 C(\gamma,s)Te^{\nu\gamma^2 T}\left(\hs{\omega_0}{s}+2T\hs{v_0}{s+1}\right)\right)\nonumber\\
&\le& 2\hs{v_0}{s+1}, \quad \forall t\in [0,T].\\
\nonumber\\
\hs{\omega^{k+1}}{s}&\le & C(\gamma,s)e^{\nu\gamma^2 T}\left(\hs{\omega_0}{s}+T\sup_{t\in [0,T]}\hs{v^{k}}{s+1}\right)\nonumber\\
&\le & C(\gamma,s)e^{\nu\gamma^2 T}\left(\hs{\omega_0}{s}+2T\hs{v_0}{s+1}\right),\quad \forall t\in [0,T].
\end{eqnarray}
Then, by induction, we prove that for any $k\ge 0$, $U^{k}=(v^k,\omega^k)$ is bounded by \myref{hs-iteration-v} and \myref{hs-iteration-w}.

We want to point that there exists $T>0$ such that the inequality \myref{T-constrain} is satisfied. One choice of $T$ is given as following:
\begin{eqnarray}
  T_1=\min\left\{\left[8C(\gamma,s)e^{\nu\gamma^2}\left(\hs{\omega_0}{s}+2\hs{v_0}{s+1}\right)\right]^{-1},1\right\}.\label{choice-T}
\end{eqnarray}

Using the choice of T in \myref{choice-T}, we can choose 
$M=2\hs{v_0}{s+1}+C(\gamma,s)e^{\nu\gamma^2}\left(\hs{\omega_0}{s}+2\hs{v_0}{s+1}\right)$, then we have
$ U^{k}\in S$, for all $k$.

Next, we will prove that $\Phi$ is a contraction mapping for some small $0<T\le T_1$.

First of all, by using Lemmas \ref{reisz} and \ref{lemma hs}, we have
\begin{eqnarray}
  \label{compress-v}
  \hs{v^{k+1}-v^{k}}{s+1}&=&\hs{\int_0^tv^k(\Bx,t')\mathcal{K}(\omega^k)_z(\Bx,t')dt'-\int_0^tv^{k-1}(\Bx,t')\mathcal{K}(\omega^{k-1})_z(\Bx,t')dt'}{s
+1}\nonumber\\
&\le& \hs{\int_0^t\left(v^k-v^{k-1}\right)(\Bx,t')\mathcal{K}(\omega^k)_z(\Bx,t')dt'}{s+1}\nonumber\\
&&+\hs{\int_0^tv^{k-1}(\Bx,t')\left(\mathcal{K}(\omega^{k})_z-\mathcal{K}(\omega^{k-1})_z\right)(\Bx,t')dt'}{s
+1}\nonumber\\
&\le & T\sup_{t\in [0, T]}\hs{v^k-v^{k-1}}{s+1}\sup_{t\in [0, T]}\hs{\mathcal{K}(\omega^k)_z}{s+1}\nonumber\\
&&+T\sup_{t\in [0, T]}\hs{v^{k-1}}{s+1}\sup_{t\in [0, T]}\hs{\mathcal{K}(\omega^{k}-\omega^{k-1})_z}{s+1}\nonumber\\
&\le & MT\left(\sup_{t\in [0, T]}\hs{v^k-v^{k-1}}{s+1}+\sup_{t\in
[0, T]}\hs{\omega^{k}-\omega^{k-1}}{s}\right).
\end{eqnarray}
On the other hand, Lemma \ref{heat-kernel} and \myref{duhamle} imply
\begin{eqnarray}
  \label{compress-omega}
  \hs{\omega^{k+1}-\omega^{k}}{s}&=&\hs{\mathcal{L}(v_z^{k},\omega_0;\Bx,t)-\mathcal{L}(v_z^{k-1},\omega_0;\Bx,t)}{s}\nonumber\\
&\le& \hs{\int_0^t\mathcal{P}(v_z^{k}-v_z^{k-1};t',t)dt'}{s}\nonumber\\
&\le &T C(\gamma,s)e^{\nu\gamma^2 T}\sup_{t\in [0,T]}\hs{v_z^{k}-v_z^{k-1}}{s} \nonumber\\
&\le & T C(\gamma,s)e^{\nu\gamma^2 T}\sup_{t\in [0,T]}\hs{v^{k}-v^{k-1}}{s+1}.
\end{eqnarray}
Let
\begin{eqnarray}
  \label{choice-T-final}
  T=\min\left\{\left[8C(\gamma,s)e^{\nu\gamma^2}\left(\hs{\omega_0}{s}+2\hs{v_0}{s+1}\right)\right]^{-1}, 
\left[2C(\gamma,s)e^{\nu\gamma^2}\right]^{-1},\;\frac{1}{2M},\;1\right\} .
\end{eqnarray}
Then, we have
\begin{eqnarray*}
\left\|U^{k+1}-U^{k}\right\|_X\le \frac{1}{2}\left\|U^{k}-U^{k-1}\right\|_X.
\end{eqnarray*}
This proves that the sequence $U^k$ converges to a fixed point of
the map $\Phi: X\rightarrow X$, and the limiting fixed point
$U=(v,\omega)$ is
a solution of the 3D model with partial viscosity. Moreover, by
passing the limit in \myref{hs-iteration-v}-\myref{hs-iteration-w},
we obtain the following {\it a priori} estimate for the solution
$v$ and $\omega$:
\begin{eqnarray}
  \label{hs-iteration-v1}
 && \hs{v}{s+1}\le 2\hs{v_0}{s+1},\\
  \label{hs-iteration-w1}
 &&\hs{\omega}{s}\le C(\gamma,s)e^{\nu\gamma^2 T}\left(\hs{\omega_0}{s}+2T\hs{v_0}{s+1}\right),
\end{eqnarray}
for $0 \leq t \leq T$ with $T$ defined in \myref{choice-T-final}.

It remains to show that the smooth solution of the 3D model
with partial viscosity is unique. Let $(v_1,\omega_1)$ and
$(v_2,\omega_2)$ be two smooth solutions of the 3D model 
with the same initial data and satisfying
$\hs{v_i}{s+1} \leq M$ and 
$\hs{\omega_i}{s} \leq M$ for $i=1,2$ and $0 \leq t \leq T$,
where $M$ is a positive constant depending on the initial data,
$\gamma$, $s$, and $T$.
Since $s > 3/2$, the Sobolev embedding theorem \cite{Evans98}
implies that 
\begin{eqnarray}
\label{embeding-vi}
&&\|v_i \|_{L^\infty(\Omega)} \leq \hs{v_i}{s+1} \leq M ,\quad i=1,2,\\
\label{embeding-wi}
&&\|{\cal{K}}(\omega_i)_z \|_{L^\infty(\Omega)} 
\leq \hs{{\cal{K}}(\omega_i)_z}{s} \leq C_s 
\hs{\omega_i}{s} \leq C_s M ,\quad i=1,2.
\end{eqnarray}

Let $v=v_1-v_2$ and $\omega = \omega_1 - \omega_2$.
Then $(v,\omega)$ satisfies
\begin{eqnarray}
&&\left\{\begin{array}{rcl}
v_t &=&4v\mathcal{K}(\omega_1)_z+4 v_2 \mathcal{K}(\omega)_z \\
\omega_t&=&v_z+\nu \Delta \omega\\
\end{array}\right., \quad (\Bx,z)\in \Omega=\Ox\times(0,\infty),
\end{eqnarray}
with 
$\omega|_{\p\Omega\backslash\Gamma}=0,\quad \left(\omega_z+\gamma \omega\right)|_\Gamma=0$, and
$\omega|_{t=0}=0$, $ v|_{t=0}=0.$
By using \myref{embeding-vi}-\myref{embeding-wi}, and 
proceeding as the uniqueness estimate for the heat equation
in \myref{unique-heat}, we can derive the following
estimate for $v$ and $\omega$: 
\begin{eqnarray}
\label{est-v2}
\frac{d}{dt} \|v\|_{L^2(\Omega)}^2 \leq C_1 (\|v\|_{L^2(\Omega)}^2
+\|\omega\|_{L^2(\Omega)}^2) ,\\
\label{est-w2}
\frac{d}{dt} \|\omega\|^2_{L^2(\Omega)} \leq C_3 (\|v\|^2_{L^2(\Omega)}
+\|\omega\|^2_{L^2(\Omega)}) ,
\end{eqnarray}
where $C_i$ ($i=1,2,3$) are positive constants depending on $M$, 
$\nu$, $\gamma$, $C_s$. In obtaining the estimate for \myref{est-w2},
we have performed integration by parts in the estimate of the
$v_z$-term in the $\omega$-equation and absorbing the contribution
from $\omega_z$ by the diffusion term. There is no contribution
from the boundary term since $v|_{z=0} = 0$.
We have also used the property 
$\|{\cal K}(\omega)_z\|_{L^2(\Omega)} \le C_s 
\|\omega\|_{L^2(\Omega)}$, which can be proved directly
by following the argument in the Appendix of \cite{HSW11}. 
Since $v_0 = 0$ and $\omega_0=0$, the Gronwall inequality 
implies that $\|v\|_{L^2(\Omega)} = \|\omega\|_{L^2(\Omega)} = 0$
for $0 \leq t \leq T$.
Furthermore, since $v \in H^{s+1}$ and $\omega \in H^s$ with
$s>3/2$, $v$ and $\omega$ are continuous. Thus
we must have $v=\omega = 0$ for $0 \leq t \leq T$. 
This proves the uniqueness of the smooth solution for the 3D model.
\end{proof}

\end{document}